\documentclass[leqno,11pt]{amsart}
\usepackage{amssymb}
\usepackage{amssymb, amsmath,txfonts}
\usepackage{tikz}
\usepackage{multirow}

\let\al=\alpha

\let\d=\delta
\let\e=\varepsilon

\let\f=\frac

\let\D=\Delta

\let\tri=\triangle

\def\cD{{\mathcal D}}

\def\cL{{\mathcal L}}
\def\cM{{\mathcal M}}

\def\cP{{\mathcal P}}

\def\C{\mathop{\mathbb C\kern 0pt}\nolimits}
\def\DD{\mathop{\mathbb D\kern 0pt}\nolimits}
\def\EE{\mathop{{\mathbb E \kern 0pt}}\nolimits}
\def\K{\mathop{\mathbb K\kern 0pt}\nolimits}
\def\N{\mathop{\mathbb N\kern 0pt}\nolimits}
\def\Q{\mathop{\mathbb Q\kern 0pt}\nolimits}
\def\R{\mathop{\mathbb R\kern 0pt}\nolimits}
\def\SS{\mathop{\mathbb S\kern 0pt}\nolimits}
\def\ZZ{\mathop{\mathbb Z\kern 0pt}\nolimits}
\def\TT{\mathop{\mathbb T\kern 0pt}\nolimits}

\def\bh{\overline{h}}

\def\bc{\overline{c}}
\def\bp{\overline{p}}
\def\bb{\overline{b}}

\def\NLS{nonlinear Schr\"odinger equation\ }
\def\iff{if and only if}

\def\p{\partial}
\def\th{\theta}

\newcommand{\beq}{\begin{equation}}
\newcommand{\eeq}{\end{equation}}
\newcommand{\ben}{\begin{eqnarray}}
\newcommand{\een}{\end{eqnarray}}
\newcommand{\beno}{\begin{eqnarray*}}
\newcommand{\eeno}{\end{eqnarray*}}

\newtheorem{defi}{Definition}[section]
\newtheorem{thm}{Theorem}[section]

\newtheorem{rmk}{Remark}[section]
\newtheorem{col}{Corollary}[section]
\newtheorem{prop}{Proposition}[section]

\def\beginpf {\noindent {\it Proof. }}
\def\endpf {\hfill $\Box$}

\begin{document}

\title{Large time blow up for a perturbation of the cubic Szeg\H{o} equation}
\author[H. XU]{Haiyan XU}
\address{Haiyan XU\\
Laboratoire de Math\'ematiques d'Orsay\\
Universit\'e Paris-Sud (XI)\\
91405, Orsay Cedex, France}

\email{haiyan.xu@math.u-psud.fr}

\subjclass[2010]{ 35B44, 37J35, 47B35}

\begin{abstract}
We consider  the following Hamiltonian equation on a special manifold of rational functions,
\[i\p_tu=\Pi(|u|^2u)+\al (u|1),\ \al\in\R,\]
where $\Pi $ denotes the Szeg\H{o} projector on the Hardy space of the circle $\SS^1$. The equation with $\al=0$ was first introduced by G\'erard and Grellier in \cite{GG1} as a toy model for totally non dispersive evolution equations. We establish the following properties for this equation. For $\al<0$, any compact subset of initial data leads to a relatively compact subset of trajectories. For $\al>0$, there exist trajectories on which high Sobolev norms exponentially grow in time.
\end{abstract}

\keywords{Szeg\H{o} equation; integrable Hamiltonian systems; Lax pair; large time blow up}
\footnote{This work was supported by grants from R\'egion Ile-de-France.}
\maketitle

\section{Introduction}
The study on the long time behavior of solutions of Schr\"odinger type Hamiltonian equations is a central issue in the theory of dispersive nonlinear partial differential equations. For instance, Colliander, Keel, Staffilani, Takaoka and Tao studied the following cubic defocusing \NLS in \cite{CKSTT},
\beq\label{NLS} i\p_t u+\tri u=\pm |u|^2u\ , \quad (t,x)\in\R\times \TT^2\ .\eeq
In that paper, they constructed solutions with small $H^s$ norm at the initial moment, which present a large Sobolev $H^s$ norm at a sufficiently long time $T$. Guardia and Kaloshin improved this result by refining the estimates on the time $T$ \cite{GK}. Zaher Hani studied a version of \NLS obtained by canceling the least resonant part, and showed the existence of unbounded trajectories in high Sobolev norms \cite{ZH}. Recently, Hani, Pausader, Tzvetkov and Visciglia studied the \NLS (\ref{NLS}) on the spatial domain $\R\times\TT^d$, and obtained global solutions to the defocusing and focusing problems on (for any $d\ge 2$) with infinitely growing high Sobolev norms $H^s$ \cite{HPTV}.

There is another related result by G\'erard and Grellier \cite{GG3}. They considered the following degenerate half wave equation on the one dimensional torus,
\beq\label{hw} i\p_t u-|D|u=|u|^2u\ .\eeq
They found solutions with small Sobolev norms at initial time which become much larger as time grows. More precisely, there exist sequences of solutions $u^n$ and $t^n$ such that $\|u_0^n\|_{H^r}\to0$ for any $r$, but
\[\|u^n(t^n)\|_{H^s}\sim\|u_0^n\|_{H^s}\left(\log\f{1}{\|u_0^n\|_{H^s}}\right)^{2s-1},\ s>1\ .\]
In fact, the above result is a consequence of the studies on the so-called {\em cubic Szeg\H{o} equation} which is introduced by G\'erard and Grellier as a model of non-dispersive dynamics \cite{GG1,GG2},
\beq\label{szego} i\p_t u=\Pi(|u|^2u)\ .\eeq
The above equation turns out to be the resonant part of the half wave equation (\ref{hw}). The operator $\Pi$, which is the so-called Szeg\H{o} operator, is defined as a projector onto the non-negative frequencies. If $u\in\cD'(\SS^1)$ is a distribution on the circle $\SS^1=\{z\in \C:\  |z|=1\}$, then
\beq\label{proj}\Pi(u)=\Pi\big(\sum\limits_{k\in\ZZ}\hat{u}(k)e^{ik\th}\big)=\sum\limits_{k\geq0}\hat{u}(k)e^{ik\th}.\eeq
Notice that, on the Hilbert space $L^2(\SS^1)$ endowed with the inner product
\beq\label{innprod}(u\ |\ v)=\f1{2\pi}\int_{-\pi}^{\pi}u(e^{ix})\overline{v(e^{ix})}dx\ ,\eeq
$\Pi$ is the orthogonal projector on the subspace $L^2_+(\SS^1)$ defined by the conditions
\[\forall k<0,\ \hat{u}(k)=0\ .\]

G\'erard and Grellier studied the Szeg\H{o} equation on the space $H^{\f12}(\SS^1)\cap L^2_+(\SS^1):=H^{\f12}_+(\SS^1)$ and displayed two Lax pair structures for this completely integrable system \cite{GG1,GG2}. Moreover, they established an explicit formula of every solution with rational initial data \cite{GG4} and illustrated the large time behavior of Sobolev norms of the solutions, for instance,
\begin{thm}\label{szgg}\cite{GG1}
Every solution $u$ of (\ref{szego}) on
$$\widetilde{\cM}(1):=\left\{u=\f{a+bz}{1-pz}:\ 0\ne a\in\C, b \in\C, p\in \C, |p| < 1 , a + bp \ne 0\right\}$$
satisfies
\[\forall s>\f12,\ \sup\limits_{t\in\R}\|u(t)\|_{H^s}<\infty.\]
However, there exists a family of Cauchy data $u_0^\e$ in $\widetilde{\cM}(1)$ which converges in $\widetilde{\cM}(1)$ for the $C^\infty(\SS^1)$ topology as $\e\to0$, and $K>0$ such that
the corresponding solutions of (\ref{szego}) $u^\e$ satisfy
\[\forall\e>0,\ \exists t^\e>0,\ \|u^\e(t^\e)\|_{H^s}\ge K(t^\e)^{2s-1} \text{ as } t^\e\to\infty,\ \forall s >\f12.\]
\end{thm}

Another result on this Szeg\H{o} equation was obtained by Pocovnicu \cite{PO1,PO2}, who studied this equation by replacing the circle $\SS^1$ with the real line and got a polynomial growth of high Sobolev norms (Corollary 4, \cite{PO2}), which says that there exists a solution $u$ of the Szeg\H{o} equation and a constant $C>0$ such that $\|u(t)\|_{H^s}\ge C|t|^{2s-1}$ for sufficiently large $|t|$.

The aim of this manuscript is to study the properties of global solutions for the following Hamiltonian equation on $L^2_+(\SS^1)$, which is the cubic Szeg\H{o} equation with a linear perturbation,
\beq\label{alsz}\left\{\begin{split}
&i\p_t u=\Pi(|u|^2u)+\al(u\ |\ 1),\ \al\in\R\ ,\\
&u(0,x)=u_0(x)\ ,
\end{split}\right.\eeq
Recall that, in view of the above definition (\ref{innprod}),
\[(u\ |\ 1)=\f1{2\pi}\int_{-\pi}^{\pi}u(e^{ix})dx\]
is the average of $u$ on $\SS^1$.

The equation (\ref{alsz}), called the $\al-$Szeg\H{o} equation, inherits three formal conservation laws:
\[\begin{aligned}
&\text{mass:  }Q(u):=\int_{\SS^1}|u|^2 \f{d\th}{2\pi}=\|u\|_{L^2}^2\ ,\\
&\text{momentum:  }M(u):=(Du\ |\ u),\ D:=-i\p_\th=z\p_z\ ,\\
&\text{energy:  }E_\al(u):=\f{1}{4}\int_{\SS^1}|u|^4\f{d\th}{2\pi}+\f12\al|(u|1)|^2\ .
\end{aligned}\]

Slight modifications of the proof of the well-posedness result in \cite{GG1} lead to the result that the $\al-$Szeg\H{o} equation is globally well-posed in $H^s_+(\SS^1)=H^s(\SS^1)\cap L^2_+(\SS^1)$ for $s\ge\f12$ as follows:
\begin{thm}\label{wp}
Given $u_0\in H^{\f12}_+(\SS^1)$, there exists a unique global solution $u\in C(\R;H^{\f12}_+)$ of (\ref{alsz}) with $u_0$ as the initial condition. Moreover, if $u_0\in H^s_+(\SS^1)$ for some $s>\f12$, then $u\in C^\infty(\R;H^s_+)$. Furthermore, if $u_0\in H^s_+(\SS^1)$ with $s>1$, the Wiener norm of $u$ is bounded uniformly in time,
\beq\sup\limits_{t\in\R}\|u(t)\|_W:=\sup\limits_{t\in\R}\sum\limits_{k=0}^\infty |\widehat{u(t)}(k)|\le C_s\|u_0\|_{H^s}\ .\eeq
\end{thm}

Now, we present our main results. In our case with a perturbation term, we gain the following statement that for the case $\al<0$ the Sobolev norm stays bounded uniformly in time, while for $\al>0$, it may grow exponentially fast.

\begin{thm}\label{mainthm}
Let $u_0=b_0+\f{c_0 z}{1-p_0z}$, $c_0\ne0$, $|p_0|<1$.

For $\al<0$, the Sobolev norm of the solution will stay bounded,
\beq\|u(t)\|_{H^s}\le C,\text{ $C$ does not depend on time $t$, }s\ge0\ .\eeq

For $\al>0$, the solution $u$ of the $\al-$Szeg\H{o} equation (\ref{alsz}) has a Sobolev norm growing exponentially in time,
\beq\label{mainest}\|u(t)\|_{H^s}\simeq e^{C_{\al,s}|t|},\ s>\f12,\ C_{\al,s}>0,\ |t| \rightarrow \infty\ ,\eeq
\iff\
\beq E_\al=\f{1}{4}Q^2+\f{\al}{2}Q\ .\eeq

\end{thm}

\begin{rmk}\label{firstrmk}

Here are several remarks:

\begin{enumerate}
  \item Together with the results in \cite{GG1,GG2}, we now have a complete picture for the high Sobolev norm of the solutions to the $\alpha$-Szeg\H{o} equation. For $\al<0$, it stays bounded (uniformly on time), for $\al>0$, it turns out to have an exponential growth for some initial data satisfying the condition in the Theorem \ref{mainthm}. Finally, for $\al=0$, the trajectories of the Szeg\H{o} equation with rational initial data are quasiperiodic with instability of the $H^s$ norm as in Theorem \ref{szgg}.
  \item Our result is in strong contrast with Bourgain's and Staffilani's results for the dispersive equations in \cite{Bg2,Sta}, which say that the dispersive equations admit polynomial upper bounds on Sobolev norms. Here, we give an example of exponential growth of Sobolev norms for a non dispersive model.
  \item The solutions to the $\al$-Szeg\H{o} equation admit an exponential upper bound of the Sobolev norms. Assume $s>1$, it is easy to solve (\ref{alsz}) locally in time. More precisely, one has to solve the integral equation
      \[u(t)=u_0-i\int_0^t \Big(\Pi(|u|^2u)+\al(u|1)\Big)dt'\ .\]
      Thus
      \[\|u(t)\|_{H^s}\le \|u_0\|_{H^s}+c\int_0^t \Big(1+\|u(t')\|_W^2\Big)\|u(t')\|_{H^s}dt'\ ,\]
      since by Theorem \ref{wp}, the Wiener norm is uniformly bounded, then by Gronwall's inequality, we have
      \[\|u(t)\|_{H^s}\le \|u_0\|_{H^s}e^{ct}\ .\]
      This shows that estimate (\ref{mainest}) is the worst that can happen.
\end{enumerate}
\end{rmk}

This paper is organized as follows. In section 2, we prove that there exists a Lax pair for the $\al-$Szeg\H{o} equation based on Hankel operators. Then we define the manifolds $\cL(k):=\Big\{u:\ rk(K_u)=k,\ k\in\ZZ^+\Big\}$ with the shifted Hankel operator $K_u$. These manifolds are proved to be invariant by the flow and can be represented as sets of rational functions. In this paper we will just consider the solutions $u\in \cL(1)$. We plan to address the other cases in a forthcoming work. In section 3, we prove the large time blow up result and the boundedness of the Wiener norm to show that our result is optimal. Furthermore, we provide an example which describes the energy cascade. Finally, we present some perspectives in section 4.

\section{The Lax pair structure}
For $u\in E\subset\cD'(\SS^1)$, we define $E_+$ by canceling the negative Fourier modes of $u$,
\[E_+=\Big\{u\in E:\  \forall k<0, \hat{u}(k)=0\Big\}\ .\]
In particular, $L_+^2$ is the Hardy space of $L^2$ functions which extend to the unit disc $D=\big\{z\in\C:\ |z|<1\big\}$ as holomorphic functions
\[ u(z)=\sum\limits_{k\geq0}\hat{u}(k)z^k,\quad\sum\limits_{k\geq0}|\hat{u}(k)|^2<\infty\ .\]
An element of $L^2_+$ can therefore be seen either as a square integrable  function $u=u(e^{i\theta})$ on the circle with only nonnegative Fourier modes, or a holomorphic function $u=u(z)$ on the
unit disc with  square summable Taylor coefficients.

Using the Szeg\H{o} projector defined as (\ref{proj}), we first introduce two important classes of operators on $L^2_+(\SS^1)$, namely, the Hankel and Toeplitz operators.

Given $u\in H^{\f12}_+(\SS^1)$, a Hankel operator $H_u:\ L^2_+\to L^2_+$ is defined by
\[H_u(h)=\Pi(u\bh)\ .\]
Notice that $H_u$ is $\C-$antilinear and symmetric with respect to the real scalar product $\mathrm{Re}(u|v)$. In fact, it satisfies
\[(H_u(h_1)\ |\ h_2)=(H_u(h_2)\ |\ h_1)\ .\]
Moreover, $H_u$ is a Hilbert-Schmidt operator with
\[\mathrm{Tr}(H_u^2)=\sum\limits_{n=0}^\infty(n+1)|\hat{u}(n)|^2\ .\]

Given $b\in L^\infty(\SS^1)$, a Toeplitz operator $T_b:\ L^2_+\to L^2_+$ is defined by
\[T_b(h)=\Pi(bh)\ .\]
$T_b$ is $\C-$linear, bounded, and is self-adjoint if and only if $b$ is real-valued.

The cubic Szeg\H{o} equation was proved to admit two Lax pairs as follows:
\begin{thm}[\cite{GG1}, Theorem 3.1]
Let $u\in C(\R,H^s(\SS^1))$ for some $s>\f12$. The cubic Szeg\H{o} equation
\beq\label{szlax} i\p_tu=\Pi(|u|^2u)\eeq
has two Lax pairs $(H_u, B_u)$ and $(K_u, C_u)$, namely, if $u$ solves (\ref{szlax}), then
\beq\f{d H_u}{dt}=[B_u,H_u]\ , \ \f{d K_u}{dt}=[C_u,K_u]\ ,\eeq
where
\[B_u=\f{i}{2}H_u^2-i T_{|u|^2}\ , \ K_u:=T_z^*H_u\ , \ C_u=\f{i}{2}K_u^2-iT_{|u|^2}\ .\]
\end{thm}

\begin{col}
The perturbed Szeg\H{o} equation (\ref{alsz}) with $\al\ne0$ still has one Lax pair $(K_u, C_u)$.
\end{col}

\beginpf
The proof is based on the following identity (\cite{GG4}, Lemma 1),
\beq\label{Hpi}H_{\Pi(|u|^2u)}=T_{|u|^2}H_u+H_u T_{|u|^2}-H_u^3\ .\eeq
Using equation (\ref{alsz}) and (\ref{Hpi}),
\[\f{d H_u}{dt}=H_{-i\Pi(|u|^2u)-i\al(u|1)}=-i(T_{|u|^2}H_u+H_u T_{|u|^2}-H_u^3)-i\al(u\ |\ 1)H_1\ .\]
Using the anti-linearity of $H_u$, we deduce that
\beq\label{dH}\f{d H_u}{dt}=[B_u,H_u]-i\al(u\ |\ 1)H_1\ ,\eeq
which means that $(H_u,B_u)$ is no longer a Lax pair. Fortunately, we have $T_z^*H_1=0$, which leads to the following identity
\beq\label{dK}\f{d K_u}{dt}=[C_u,K_u]\ .\eeq
\endpf

An important consequence of this Lax pair structure is the existence of finite dimensional submanifolds of $L^2_+(\SS^1)$ which are invariant by the flow of (\ref{alsz}). To describe these manifolds, G\'erard and Grellier (Appendix 4, \cite{GG1}) proved a Kronecker-type theorem that, the Hankel operator $H_u$ is of finite rank $k$ if and only if $u$ is a rational function of the complex variable $z$, with no poles in the unit disc, and of the form $u(z)=\f{A(z)}{B(z)}$ with $A\in\C_{k-1}[z]$, $B\in\C_{k}[z]$, $B(0)=1$, $\deg(A)=k-1$ or $\deg(B)=k$, $A$ and $B$ have no common factors and $B(z)\ne0$ if $|z|\le1$. In fact, we can prove a similar theorem for our case.

\begin{defi}
Let $k$ be a positive integer, we define
\beq\cL(k):=\left\{u\in H^{\f12}_+(\SS^1):\ rk(K_u)=k\right\}\ .\eeq
\end{defi}

Due to the Lax pair structure, the manifolds $\cL(k)$ are invariant by the flow.

\begin{thm}
$u\in \cL(k)$ if and only if u is a rational function satisfying
\[u(z)=\f{A(z)}{B(z)}\text{ with }\ A,B\in\C_k[z],A\wedge B=1, \deg( A)= k \text{ or }
\deg( B)= k, B^{-1}(\{0\})\cap \overline{D}=\emptyset\ ,\]
where $A\wedge B=1$ means $A$ and $B$ have no common factors.
\end{thm}

\beginpf
The proof is based on the results by G\'erard and Grellier (see Appendix 4, \cite{GG1}), they proved that
\[\begin{split}\cM(k+1)&=\{u:\ rk(H_u)=k+1\}\\
&=\left\{\begin{split}&u(z)=\f{A(z)}{B(z)}:\ A\in\C_{k}[z], B\in\C_{k+1}[z], B(0)=1,\deg(A)=k\\
&\text{ or } \deg(B)= k+1, A\wedge B=1,\ B^{-1}(0)\cap \overline{D}=\emptyset\end{split}\right\}.\end{split}\]
For $u\in \cM(k+1)$, $\dim\mathrm{Im} H_u=k+1$, then $u,T_z^*u,\cdots,(T_z^*)^{k+1}u$ are linearly dependent,
i.e, there exist $C_\ell$, not all zero, such that $\sum\limits_{\ell=0}^{k+1}C_\ell(T_z^*)^\ell u=0$. We get\[\sum\limits_{\ell=0}^{k+1}C_\ell\hat u(\ell+n)=0\ ,\ \forall n\geq 0\ .\]
This is a recurrent equation for sequence $\hat u$. It can be solved by means of elementary linear algebra. Define
\[P(X)=\sum\limits_{\ell=0}^{k+1}C_\ell X^\ell=C\prod\limits_{p\in\cP}(X-p)^{m_p}\ ,\]
where $\cP=\{p\in\C: P(p)=0\}$ and $m_p$ is the multiplicity of $p$.\\
$(\hat u(n))_{n\geq0}$ is a linear combination of the following sequences:
\[n^\ell p^{n-\ell}, \ p\neq 0, \ 0\leq \ell\leq m_p-1\ ,\]
\[\d_{nm},\ p=0, \ 0\leq m\leq m_0-1\ .\]
Recall that\[u(z)=\sum\limits_{n\geq0} \hat u(n)z^n \text{ for }  |z|<1\ ,\]
then $u$ is a linear combination of $\f{1}{(1-pz)^{\ell+1}}$ with $0<|p|<1$ for $ 0\leq
\ell\leq m_p-1$, and of $z^\ell$ for $0\leq \ell\leq m_0-1$.\\
Consequently, $u(z)=\f{A(z)}{B(z)}$ with
 \[\begin{split}&\deg( A) \leq k,\ \deg( B) =k+1,\text{ if }p\neq0,\ p\in \cP\ ,\\
 &\deg(A) = k,\ \deg(B) \leq k,\text{ if }0\in \cP\ .\end{split}\]
Note that\[0\in\cP\] is equivalent to \[1\in \mathrm{Im} H_u\] or to \[\ker K_u \cap
\mathrm{Im} H_u \neq \{0\}\ ,\]
since $K_u=T_z^*H_u$, $rk(H_u)-1\le rk(K_u)\le rk(H_u)$. For $u\in\cL(k)$, $rk(K_u)=k$,
then $u=\f{A(z)}{B(z)}$ with
\beno\begin{split}&\deg( A)\le k-1, \deg( B)=k, \text{ if $rk (H_u)=rk (K_u)=k$}\ , \\
&\deg( A)= k, \deg (B)\le k, \text{ if $rk (H_u)=rk (K_u)+1=k+1$}\ .\end{split}\eeno
The proof of the converse is similar. So
\[\begin{split}\cL(k)&=\{u:\ rk(K_u)=k+1\}\\
&=\left\{\begin{split}&u(z)=\f{A(z)}{B(z)}:\ A\in\C_{k}[z], B\in\C_{k}[z], B(0)=1,\deg(A)=k\\
&\text{ or } \deg(B)= k, A\wedge B=1,\ B^{-1}(0)\cap \overline{D}=\emptyset\end{split}\right\}\ .\end{split}\]
The proof is completed.
\endpf

\section{The proof of the main theorem}
In this section, we will prove that the $\al-$Szeg\H{o} equation (\ref{alsz}) admits the large time blow up as in Theorem \ref{mainthm}, we will also give an example to describe this phenomenon in terms of energy transfer to high frequencies. Before the proof of the main theorem, let us prove the boundedness of Wiener norm as in Theorem \ref{wp}.

\begin{prop}\label{wiener}
Assume $u_0\in H^s_+(\SS^1)$ with $s>1$, let $u$ be the corresponding unique solution of (\ref{alsz}). Then
\[\|u(t)\|_W\le C_s\|u_0\|_{H^s},\ \forall t\in\R\ .\]
\end{prop}
\beginpf
By Peller's theorem \cite{Pel}, the regularity of $u$ ensures that $H_u$ is trace class and the trace norm of $H_u$ is equivalent to the $B_{1,1}^1$ norm of $u$. Recall the definition of $B^s_{p,q}(\SS ^1)$.

Let $\chi \in C^{\infty}(\R^+)$ satisfy $\chi|_{t<1}(t)= 1$, $\chi|_{t>2}(t)= 0$, $0\leq \chi\leq 1$. Set $\psi$ as $\psi_0(t)= 1-\chi(t)$, $\psi_j(t)=\chi(2^{-j+1}t)-\chi(2^{-j}t)$. Define the operator $\D_j$ for $f\in  \cD'(\SS ^1)$ as
\[\D_jf=\sum _{k\in \ZZ}\psi _j(k)\hat f(k)\, e^{ik\th }\ .\]
Then the Besov space is defined as
\[B^s_{p,q}(\SS ^1):=\Big\{u\in\cD'(\SS ^1):2^{js}\|\D_jf\|_{L^p}\in l_j^q,1\leq p, q\leq +\infty,\ 0\le j\le +\infty\Big\}\ ,\]
with the norm $\|u\|_{B^s_{p,q}(\SS ^1)}=\left(\sum\limits_{j=0}^{+\infty} (2^{js}\|\D_jf\|_{L^p})^q\right)^{\f{1}{q}}$.

Observe that there exist $C$, $C_s\ >0$, such that
\beq\begin{split}\|u\|_{B^1_{1,1}}&=\sum\limits_{j=0}^{+\infty} 2^{j}\|\D_ju\|_{L^1}\le C\sum\limits_{j=0}^{+\infty} 2^{j}\|\D_ju\|_{L^2}\\
&\le C\left(\sum\limits_{j=0}^{+\infty} 2^{2js}\|\D_ju\|_{L^2}^2\right)^\f12\left(\sum\limits_{j=0}^{+\infty} 2^{2j(1-s)}\right)^\f12\\
&\le C_s\|u\|_{H^s},\ \forall s>1\ .\end{split}\eeq

So for $u\in H^s$ with $s>1$, $H_u$ is trace class, and
\[\mathrm{Tr}(|H_u|)\le C_s\|u\|_{H^s}\ .\]
Since $K_u=T_z^*H_u$, then
\[K_u^2=H_u^2-(\cdot\ |\ u)u\ ,\]
then
\[\mathrm{Tr}(|K_u|)\le \mathrm{Tr}(|H_u|)\ .\]
Due to the Lax pair structure, we get $K_{u(t)}$ is isospectral to $K_{u_0}$, then
\[\mathrm{Tr}(|K_{u(t)}|)=\mathrm{Tr}(|K_{u_0}|)\ ,\]
so
\[\mathrm{Tr}(|K_{u(t)}|)\le C_s\|u_0\|_{H^s}\ .\]
Since $\|u\|_W=|\hat{u}(0)|+\sum\limits_{n\ge1}|\hat{u}(n)|$ and $|\hat{u}(0)|\le\|u\|_{L^2}$, we just need to show that
\[\sum\limits_{n\ge1}|\hat{u}(n)|\le C \mathrm{Tr}(|K_u|)\ .\]

Let $e_n$ as the orthonormal basis of $L_+^2$, then for any bounded operator $B$,
\[\sum\limits_n\Big|(K_u e_n\ |\ Be_n)\Big|\le \mathrm{Tr}(|K_u|)\|B\|.\]
Then we gain that $\sum\limits_{n\ge1}|\hat{u}(2n)|+\sum\limits_{n\ge1}|\hat{u}(2n+1)|\le \mathrm{Tr}(|K_u|)$, by taking $B=T_z$ and $B=\mathrm{Id}$. This completes the proof.
\endpf

\begin{rmk}
In fact, to prove the global wellposedness, it is natural to use the Brezis-Gallou\"et type estimate (Appendix 2, \cite{GG1}), for $s>\f12$
\[\|u\|_W\le C_s\|u\|_{H^{\f12}}\left[\log(1+\f{\|u\|_{H^s}}{\|u\|_{H^{1/2}}})\right]^{\f12}\ ,\]
which leads to a double exponential on time growth for the Sobolev norm of $u$. Fortunately, by the estimate in Proposition \ref{wiener}, we know the $H^s$ norm of the solutions will admit an exponential on time upper bound for $s>1$ (see Remark \ref{firstrmk}).
\end{rmk}

Now, let us start the large time blow up theorem.
\begin{thm}\label{necsuf}
For $\al>0$, we consider the solution of the Szeg\H{o} equation (\ref{alsz}) with initial data $u_0\in\cL(1)$.
\begin{enumerate}
  \item If the trajectory issued from $u_0$ is not relatively compact in $\cL(1)$, then
    \beq\label{cond} \left|b+\f{\bp c}{1-|p|^2}\right|=\sqrt{\al}\ ,\eeq
    or equivalently
    \beq\label{inv} E_{\al}=\f{1}{4}Q^2+\f{\al}{2}Q\ .\eeq
  \item If (\ref{cond}) holds, then
    \beq\|u(t)\|_{H^s}\simeq e^{C_{\al,s}|t|},\ s>\f12,\ C_{\al,s}>0,\ |t|\to\infty\ .\eeq
\end{enumerate}
\end{thm}

\begin{rmk}
From the theorem, the equality (\ref{inv}), which is invariant by the flow, is a necessary and sufficient condition to cause the large time blow up.
\end{rmk}

\beginpf
First, since the trajectory of the solution is not relatively compact in $\cL(1)$, the level set $L(u_0):=\left\{u\in\cL(1):\ Q(u)=Q(u_0),\ M(u)=M(u_0),\ E_\al(u)=E_\al(u_0)\right\}$ is not compact in $\cL(1)$.

We rewrite $u\in\cL(1)$ as
\[u=b+\f{cz}{1-p z}\ ,\]
then the conservation laws under the coordinates $b,\ p,\ c$ are given as
\[\begin{aligned}
&Q=\|u\|_{L^2}^2=\f{|c|^2}{1-|p|^2}+|b|^2\ ,\\
&M=(Du\ |\ u)=\f{|c|^2}{(1-|p|^2)^2}\ ,\\
&E_\al=\f{1}{4}\|u\|_{L^4}^4+\f{\al}{2}|(u|1)|^2\\
&=\f{1}{4}\left[|b|^4+\f{4|b|^2|c|^2}{1-|p|^2}+\f{|c|^4(1+|p|^2)}{(1-|p|^2)^3}
+\f{4|c|^2\mathrm{Re}(bp\bc)}{(1-|p|^2)^2}\right]+\f{\al}{2}|b|^2\ .
\end{aligned}\]
If $u\in\cL(1)$ stays in a compact of $\cL(1)$ \iff\ $|b|\le C$, $\f1C\le|c|\le C$ and $|p|\le k<1$ with some constant $C$ and $k$. Otherwise, due to the formulas of mass $Q$ and momentum $M$, there exist $t_n\to\infty$ such that $|c(t_n)|$ and $1-|p(t_n)|^2$ tend to 0 at the same order. Using the formula of $Q$ and $E_\al$, we have
\[|b(t_n)|^2\to Q,\ \f{1}{4}|b(t_n)|^4+\f{\al}{2}|b(t_n)|^2\to \ E_\al\ .\]

Since the limit should be unique,
\[ E_\al=\f{1}{4}Q^2+\f{\al}{2}Q\ .\]

Using the formula of mass and energy, (\ref{inv}) can be rewritten under coordinates of $b$, $p$, $c$ as
\[|b|^2+\f{|c|^2|p|^2}{(1-|p|^2)^2}+2\mathrm{Re} \left(\f{bp\bc}{1-|p|^2}\right)=\al\ ,\]
simplify the left hand side, we get
\[\Big|b+\f{\bp c}{1-|p|^2}\Big|=\sqrt{\al}\ .\]

Now, we turn to prove that (\ref{cond}) is sufficient to cause the exponential growth of Sobolev norms.
Writing as before
\[u(t)=b(t)+\f{c(t)z}{1-p(t)z}\ ,\]
then the terms $\p_t u$, $\Pi(|u|^2u)$, $(u|1)$ can be represented as linear combinations of $1$, $\f{z}{1-pz}$ and $\f{z^2}{(1-pz)^2}$,
\[\left\{
\begin{split}
\p_t u&= \p_t b+\p_t c\f{z}{1-pz}+\p_t p\f{z^2}{(1-pz)^2}\ ,\\
\Pi(|u|&^2u)= |b|^2b+\f{2b|c|^2}{1-|p|^2}+\f{|c|^2c\bp}{1-|p|^2}\\
&+ \left[2|b|^2c+\f{2b|c|^2p}{1-|p|^2}+\f{1+|p|^2}{1-|p|^2}|c|^2c\right]\f{z}{1-pz}\\
&+ \left[c^2\bb+\f{|c|^2cp}{1-|p|^2}\right]\f{z^2}{(1-pz)^2}\ ,\\
(u\ |\ &1)=b\ .
\end{split}\right.\]
then (\ref{alsz}) reads
\beq\label{ode}
\left\{
\begin{aligned}
&i\p_t b = |b|^2b+\f{2b|c|^2}{1-|p|^2}+\f{|c|^2c\bp}{(1-|p|^2)^2}+\al b\ ,\\
&i\p_t c = 2|b|^2c+\f{2b|c|^2p}{1-|p|^2}+\f{|c|^2c}{(1-|p|^2)^2}\ ,\\
&i\p_t p = c\bb+\f{|c|^2p}{1-|p|^2}\ .
\end{aligned}
\right.
\eeq

Using the second equation of (\ref{ode}), we gain
\beq\label{dc1} \f{d|c|}{dt}=\f{2|c|}{1-|p|^2}\mathrm{Im} (bp\bc)\ ,\eeq

The equality together with (\ref{cond}) gives us
\[\begin{split}
\left(\f{d|c|}{|c|dt}\right)^2&=
\f{4(\mathrm{Im} (bp\bc))^2}{(1-|p|^2)^2}\\
&=\f{4|bp\bc|^2}{(1-|p|^2)^2}-\f{4(\mathrm{Re} (bp\bc))^2}{(1-|p|^2)^2}\\
&=\f{4|bp\bc|^2}{(1-|p|^2)^2}-\left[\al-|b|^2-\f{|c|^2|p|^2}{(1-|p|^2)^2}\right]^2\\
&=\f{4|bp\bc|^2}{(1-|p|^2)^2}-\left[\al-|b|^2-\f{|c|^2}{(1-|p|^2)^2}+\f{|c|^2}{1-|p|^2}\right]^2\\
&=\f{4|bp\bc|^2}{(1-|p|^2)^2}-\left[\al-Q-M+2\f{|c|^2}{1-|p|^2}\right]^2\\
&=\f{4|bp\bc|^2}{(1-|p|^2)^2}-\f{4|c|^4}{(1-|p|^2)^2}-\f{4|c|^2}{1-|p|^2}\left[\al-|b|^2-\f{|c|^2}
{(1-|p|^2)^2}-\f{|c|^2}{1-|p|^2}\right]-(\al-Q-M)^2\\
&=\f{4|b|^2|c|^2}{(1-|p|^2)^2}+\f{4|c|^4}{(1-|p|^2)^3}-\al\f{4|c|^2}{1-|p|^2}-(\al-Q-M)^2\\
&=4\left(|b|^2+\f{|c|^2}{1-|p|^2}\right)\f{|c|^2}{(1-|p|^2)^2}-\al\f{4|c|^2}{1-|p|^2}-(\al-Q-M)^2\\
&=4QM-4\al\sqrt{M}|c|-(\al-Q-M)^2\ .
\end{split}\]
Thus
\[\left(\f{d\log|c|}{dt}\right)^2=-4\al\sqrt{M}|c|+4QM-(\al-M-Q)^2\ .\]

Since $0\le|c|\le1$, then $c_{\al,M,Q}\le (\f{d|c|}{|c|dt})^2\le C_{\al,M,Q}$, which leads to an exponential decay in time for $|c|$,
\[|c|(t)\simeq|c(0)|e^{-C|t|}\ ,\]
with the positive constant $C$ depending on $\al$ and $M$, $Q$.

Notice that $\hat{u}(k,t)=cp^{k-1}$ for $k\ge1$, using Fourier expansion, we obtain, as $|p|$ approaches 1,
\[\|u\|_{H^s}^2\simeq\f{|c|^2}{(1-|p|^2)^{2s+1}}\ .\]
Since $M(u)=\f{|c|^2}{(1-|p|^2)^2}=$ constant, we get
$\|u\|_{H^s}^2\simeq |c|^{-(2s-1)}\simeq e^{C(2s-1)|t|}$, which has an exponential growth as $s>\f12$. The proof is complete.

\endpf

\begin{col}
We do not have the growth of $H^s$ norms for small data in $\cL(1)$. In other words, if $\|u(0)\|_{H^{\f12}_+}<<\sqrt{\al}$, the higher Sobolev norm will never grow to infinity.
\end{col}

\beginpf
$\|u(0)\|_{H^{\f12}_+}<<\sqrt{\al}$, then
\[\left|b+\f{c\bp}{1-|p|^2}\right|\le\sqrt{Q}+\sqrt{M}\lesssim\|u(0)\|_{H^{\f12}_+}<<\sqrt{\al}\ .\]
According to the necessary and sufficient condition (\ref{cond}), there is no norm explosion.
\endpf

\begin{rmk}
Consider a family of Cauchy data given by
\[u_0^\e=z+\e,\ \e\in\C \text{ and } \e\ne\sqrt{\al}\ .\]

For the case $\al=0$, G\'erard and Grellier got the following instability of $H^s$ norms
\[\|u^\e(t^\e)\|_{H^s}\simeq\e^{-(2s-1)},\ s>\f12\ .\]
However, we do not have such an instability result for $\al>0$. In fact, using the theorem \ref{necsuf}, we know there exists a constant $C=C(\al)$ such that,
\[\sup\limits_{\e\ne\sqrt{\al}}\sup\limits_{t\in\R}\|u^\e(t)\|_{H^s}<C\ .\]
\end{rmk}

Now, we give an example to display the energy cascade in Theorem \ref{necsuf}.
\begin{thm}
Given $\al>0$.
\beq\label{sz}
\left\{
\begin{aligned}
&i\p_t u = \Pi(|u|^2 u) + \al (u\ |\ 1)\ ,\\
&u|_{t=0}=z+\sqrt{\al}, \quad z\in\SS^1\ .
\end{aligned}
\right.
\eeq
For all $s>\f12$, the above equation is globally well-posed in $H^s$ and the solution satisfies
\[\|u(t)\|_{H^s}\simeq e^{(2s-1)\sqrt{\al}t},\ t\to\infty\ .\]
\end{thm}

\beginpf
Firstly, since $u_0=z+\sqrt{\al}$, the conserved quantities are
\[Q=1+\al,\ M=1,\ E_\al=\f{1}{4}(1+\al)(1+3\al)\ ,\]
then $u_0\in\cL(1)$.
So by the proof of Theorem \ref{necsuf},
\[\left(\f{d}{dt}|c|\right)^2=4\al|c|^2(1-|c|)\ .\]
Together with the initial condition $|c|(0) = 1$, we get for $t>0$ (same strategy for $t<0$),
\beq\label{dc}\f{d}{dt}|c|=-2\sqrt{\al}|c|\sqrt{1-|c|}\ .\eeq
\[|c|(t)=\f{4e^{2\sqrt{\al t}}}{(1+e^{2\sqrt{\al t}})^2}\ .\]
By (\ref{cond}), we can get
\[\mathrm{Re}(bp\bc)=|c|^2-|c|\ ,\]
and by (\ref{dc1}) and (\ref{dc}), we have
\[\mathrm{Im}(bp\bc)=-\sqrt{\al}|c|\sqrt{1-|c|}\ ,\]
so
\[bp\bc=\mathrm{Re}(bp\bc)+i\mathrm{Im}(bp\bc)=|c|^2-|c|-i\sqrt{\al}|c|\sqrt{1-|c|}\ .\]

The second equation of (\ref{ode}) can be simplified as follows,
\[\left\{
\begin{aligned}
&i\p_t c=\Big(1+2\al-2i\sqrt{\al}\sqrt{1-|c|}\Big)c\ ,\\
&c(0) = 1\ .
\end{aligned}
\right.\]
Then
\beq c(t)=\f{4e^{2\sqrt{\al t}}}{(1+e^{2\sqrt{\al t}})^2}e^{-i(1+2\al)t}\ .\eeq

Now, we turn to calculate $b$ and $p$, in fact, we only need to calculate their angles. Let us denote
\[b=|b|e^{i\th(t)}=\sqrt{1+\al-|c|}e^{i\th(t)}\ ,\ p=|p|e^{i\sigma(t)}=\sqrt{1-|c|}e^{i\sigma(t)}\ ,\]
then using the differential equation on $p$, we get
\[\p_t \sigma |p|=|c||p|+\mathrm{Re}\left(c\bb e^{-i\sigma}\right)=|c||p|+\mathrm{Re} \left(\f{c\bb\bp}{|p|}\right)=|c||p|+\f{1}{|p|}(|c|^2-|c|)=0\ ,\]
which means \[\sigma(t)=\sigma(0)\ .\]
Since
\[\begin{split}bp&=\f{c(bp\bc)}{|c|^2}=(|c|-1-i\sqrt{\al}\sqrt{1-|c|})e^{-i(1+2\al)t}\\
&=\sqrt{(1+\al-|c|)(1-|c|)}\left(-\f{\sqrt{1-|c|}}{\sqrt{1+\al-|c|}}-i\f{\sqrt{\al}}{\sqrt{1+\al-|c|}}\right)e^{-i(1+2\al)t}\ ,\end{split}\]
\[e^{i(\th+\sigma)}=\left(-\f{\sqrt{1-|c|}}{\sqrt{1+\al-|c|}}-i\f{\sqrt{\al}}{\sqrt{1+\al-|c|}}\right)e^{-i(1+2\al)t}\ ,\]
and $e^{i\th(0)}=1$,
thus we get \[e^{i\sigma(t)}=e^{i\sigma(0)}=e^{i(\sigma(0)+\th(0))}=-i\ ,\]
then \[e^{i\th(t)}=\left(-i\f{\sqrt{1-|c|}}{\sqrt{1+\al-|c|}}+\f{\sqrt{\al}}{\sqrt{1+\al-|c|}}\right)e^{-i(1+2\al)t}\ .\]

Finally, we have
\beq\begin{split} p(t)&=-i\sqrt{1-|c|}=-i\f{e^{2\sqrt{\al}t}-1}{e^{2\sqrt{\al}t}+1}\ ,\\
 b(t) &= \left(\sqrt{\al}-i\f{e^{2\sqrt{\al}t}-1}{e^{2\sqrt{\al}t}+1}\right)e^{-i(1+2\al)t}\ .\end{split}\eeq

Now, we get the explicit formula for the solution $u(t)=b(t)+\f{c(t)z}{1-p(t)z}$,
\beq\label{explicit formula}
\left\{
\begin{aligned}
&b(t) = \left(\sqrt{\al}-i\f{e^{2\sqrt{\al}t}-1}{e^{2\sqrt{\al}t}+1}\right)e^{-i(1+2\al)t}\,\\
&c(t) = \f{4e^{2\sqrt{\al}t}}{(1+e^{2\sqrt{\al}t})^2}e^{-i(1+2\al)t}\ ,\\
&p(t) = -i\f{e^{2\sqrt{\al}t}-1}{e^{2\sqrt{\al}t}+1}\ .
\end{aligned}
\right.
\eeq

In this case, $M(u)=\f{|c|^2}{(1-|p|^2)^2}=1$ and we get for $t\to+\infty$,
\[\|u(t)\|_{H^s}^2\simeq |c|^{-(2s-1)}\simeq C e^{2(2s-1)\sqrt{\al}t}\ .\]
\endpf

\begin{rmk}
One can illustrate this instability of Sobolev norms from the viewpoint of transfer of energy to high frequencies. The Fourier coefficients for $u=b+\f{cz}{1-pz}$ are
\[\hat{u}(k)=c(t)p(t)^{k-1},\ \forall k\ge1\ .\]
Then
\[M(u)=1=\sum\limits_{k\ge1}|k||\hat{u}(k)|^2=\sum\limits_{k\ge1}|k||c(t)|^2|p(t)|^{2(k-1)}\ .\]
With (\ref{explicit formula}), we have
\[\sum\limits_{k\ge1}\left|\f{1-e^{-2\sqrt{\al}t}}{1+e^{-2\sqrt{\al}t}}
\right|^{2k}\f{16|k|}{|(1+e^{-2\sqrt{\al}t})(1-e^{-2\sqrt{\al}t})|^2}=1\ .\]
As $t\to\infty$, we get
\[\sum\limits_{k\ge1}4|k|e^{-2\sqrt{\al}t}\exp{(-4|k|e^{-2\sqrt{\al}t})}\sim\f{1}{4}\ ,\]
so the main part of the summation is on the $k$s satisfying
\[|k|\sim e^{2\sqrt{\al}t}\ .\]
So as time becomes larger, the main part of the energy concentrates on the Fourier modes as large as $e^{2\sqrt{\al}t}$.

On the other hand, from the viewpoint of the space variable, we find that as time grows to infinity, the energy will concentrate on one point. In fact, rewrite $z=e^{ix}$, then
\beno
\begin{split}\left| u(t,x)-\sqrt{\al}-i\f{1-e^{-2\sqrt{\al}t}}{1+e^{-2\sqrt{\al}t}}\right|&=\f{|c(t)|}{|1-p(t) z|}
=\f{1-|p(t)|^2}{|1-p(t) z|}\sim \f{1-|p(t)|}{|1-p(t) z|}\\
&\sim\f{1}{\sqrt{2(e^{4\sqrt{\al}t}-1)(1-\sin x)+4}}\\
&\to 0 \text{ if and only if } x\ne \f{\pi}{2}, \ t\to\infty\ .
\end{split}\eeno
Therefore, as time tends to infinity, the value of $|u|$ will concentrate on the point $i\in\SS^1$.

Moreover, this example shows that the radius of analyticity of the solution of equation (\ref{alsz}) may decay exponentially. This shows the optimality of the result in the recent work \cite{GGT}.
\end{rmk}


Now, let us turn to the case $\al<0$.
\begin{thm}
In the case $\al<0$, for any given initial data $u_0\in\cL(1)$,
let $u=\f{a z+b}{1-p z}$ be the corresponding solution of (\ref{alsz}).
Then there exist a constant $C=C(\al)$, such that
\[\forall t,\ \|u(t)\|_{H^s}<C,\ s\ge \f12\ , \]
the constant $C>0$ is uniform for $u_0$ in a compact subset of $\cL(1)$.
\end{thm}

\beginpf
We prove this theorem by contradiction. If $u(t_n)$ would leave any compact subset of $\cL(1)$, then the Theorem \ref{necsuf} would lead to (\ref{inv}), or equivalently to the following equality,
\[\|u_0\|_{L^2}^4-\|u_0\|_{L^4}^4=2\al\Big(|(u_0\ |\ 1)|^2-\|u_0\|_{L^2}^2\Big)\ .\]
Via the Cauchy-Schwarz inequality and $\al<0$, we get
\[\|u_0\|_{L^2}=\|u\|_{L^4} \text{ and }|(u_0\ |\ 1)|=\|u_0\|_{L^2}\ ,\]
then $u_0$ should be a constant, which contradicts the fact that $u_0\in\cL(1)$.
\endpf

\section{Further studies and open problems}
In this paper, we just considered the data on the 3-(complex) dimensional manifold
\[\cL(1):=\Big\{u:\ rk K_u=1\Big\}\ .\]
It is of course natural to consider the higher dimensional case, which will be probably much more complicated. Since we have also got enough conservation laws for the case $rk K_u=2$, we have a conjecture that the system stays completely integrable for $rk K_u\ge2$. It would be interesting to know how the results of this paper extend to this bigger phase space. In particular, do small data generate large time blow up of high Sobolev norms?

\end{document}